%
%
%
%


\magnification=1200
\pretolerance=500 \tolerance=1000 \brokenpenalty=5000
\hsize=12.5cm   
\vsize=19cm
\hoffset=0.4cm
\voffset=1cm
\parskip3pt plus 1pt
\parindent=0.4cm
\def\\{\hfil\break}


\font\seventeenbf=cmbx10 at 17.28pt

\font\twelvebf=cmbx10 at 12pt
\font\eightbf=cmbx8
\font\sixbf=cmbx6

\font\eighti=cmmi8
\font\sixi=cmmi6

\font\eightrm=cmr8
\font\sixrm=cmr6

\font\eightsy=cmsy8
\font\sixsy=cmsy6

\font\eightit=cmti8
\font\eighttt=cmtt8
\font\eightsl=cmsl8

\font\seventeenbsy=cmbsy10 at 17.28pt

\font\twelvebsy=cmbsy10 at 12pt
\font\tenbsy=cmbsy10
\font\eightbsy=cmbsy8
\font\sevenbsy=cmbsy7
\font\sixbsy=cmbsy6
\font\fivebsy=cmbsy5

\font\tenmsa=msam10

\font\sevenmsa=msam7
\font\fivemsa=msam5
\newfam\msafam
  \textfont\msafam=\tenmsa
  \scriptfont\msafam=\sevenmsa
  \scriptscriptfont\msafam=\fivemsa

\font\tenmsb=msbm10
\font\eightmsb=msbm8
\font\sevenmsb=msbm7
\font\fivemsb=msbm5
\newfam\msbfam
  \textfont\msbfam=\tenmsb
  \scriptfont\msbfam=\sevenmsb
  \scriptscriptfont\msbfam=\fivemsb
\def\Bbb{\fam\msbfam\tenmsb}

\font\tenCal=eusm10
\font\sevenCal=eusm7
\font\fiveCal=eusm5
\newfam\Calfam
  \textfont\Calfam=\tenCal
  \scriptfont\Calfam=\sevenCal
  \scriptscriptfont\Calfam=\fiveCal
\def\Cal{\fam\Calfam\tenCal}

\font\teneuf=eusm10
\font\teneuf=eufm10
\font\seveneuf=eufm7
\font\fiveeuf=eufm5
\newfam\euffam
  \textfont\euffam=\teneuf
  \scriptfont\euffam=\seveneuf
  \scriptscriptfont\euffam=\fiveeuf

\font\seventeenbfit=cmmib10 at 17.28pt

\font\twelvebfit=cmmib10 at 12pt
\font\tenbfit=cmmib10
\font\eightbfit=cmmib8
\font\sevenbfit=cmmib7
\font\sixbfit=cmmib6
\font\fivebfit=cmmib5
\newfam\bfitfam
  \textfont\bfitfam=\tenbfit
  \scriptfont\bfitfam=\sevenbfit
  \scriptscriptfont\bfitfam=\fivebfit


\catcode`\@=11
\def\eightpoint{%
  \textfont0=\eightrm \scriptfont0=\sixrm \scriptscriptfont0=\fiverm
  \def\rm{\fam\z@\eightrm}%
  \textfont1=\eighti \scriptfont1=\sixi \scriptscriptfont1=\fivei
  \def\oldstyle{\fam\@ne\eighti}%
  \textfont2=\eightsy \scriptfont2=\sixsy \scriptscriptfont2=\fivesy
  \textfont\itfam=\eightit
  \def\it{\fam\itfam\eightit}%
  \textfont\slfam=\eightsl
  \def\sl{\fam\slfam\eightsl}%
  \textfont\bffam=\eightbf \scriptfont\bffam=\sixbf
  \scriptscriptfont\bffam=\fivebf
  \def\bf{\fam\bffam\eightbf}%
  \textfont\ttfam=\eighttt
  \def\tt{\fam\ttfam\eighttt}%
  \textfont\msbfam=\eightmsb
  \def\Bbb{\fam\msbfam\eightmsb}%
  \abovedisplayskip=9pt plus 2pt minus 6pt
  \abovedisplayshortskip=0pt plus 2pt
  \belowdisplayskip=9pt plus 2pt minus 6pt
  \belowdisplayshortskip=5pt plus 2pt minus 3pt
  \smallskipamount=2pt plus 1pt minus 1pt
  \medskipamount=4pt plus 2pt minus 1pt
  \bigskipamount=9pt plus 3pt minus 3pt
  \normalbaselineskip=9pt
  \setbox\strutbox=\hbox{\vrule height7pt depth2pt width0pt}%
  \let\bigf@ntpc=\eightrm \let\smallf@ntpc=\sixrm
  \normalbaselines\rm}
\catcode`\@=12

\def\eightpointbf{%
 \textfont0=\eightbf   \scriptfont0=\sixbf   \scriptscriptfont0=\fivebf
 \textfont1=\eightbfit \scriptfont1=\sixbfit \scriptscriptfont1=\fivebfit
 \textfont2=\eightbsy  \scriptfont2=\sixbsy  \scriptscriptfont2=\fivebsy
 \eightbf
 \baselineskip=10pt}

\def\tenpointbf{%
 \textfont0=\tenbf   \scriptfont0=\sevenbf   \scriptscriptfont0=\fivebf
 \textfont1=\tenbfit \scriptfont1=\sevenbfit \scriptscriptfont1=\fivebfit
 \textfont2=\tenbsy  \scriptfont2=\sevenbsy  \scriptscriptfont2=\fivebsy
 \tenbf}
        
\def\twelvepointbf{%
 \textfont0=\twelvebf   \scriptfont0=\eightbf   \scriptscriptfont0=\sixbf
 \textfont1=\twelvebfit \scriptfont1=\eightbfit \scriptscriptfont1=\sixbfit
 \textfont2=\twelvebsy  \scriptfont2=\eightbsy  \scriptscriptfont2=\sixbsy
 \twelvebf
 \baselineskip=14.4pt}

\def\seventeenpointbf{%
 \textfont0=\seventeenbf  \scriptfont0=\twelvebf  \scriptscriptfont0=\eightbf
 \textfont1=\seventeenbfit\scriptfont1=\twelvebfit\scriptscriptfont1=\eightbfit
 \textfont2=\seventeenbsy \scriptfont2=\twelvebsy \scriptscriptfont2=\eightbsy
 \seventeenbf
 \baselineskip=20.736pt}
 
       
\newbox\titlebox   \setbox\titlebox\hbox{\hfil}
\newbox\sectionbox \setbox\sectionbox\hbox{\hfil}
\def\folio{\ifnum\pageno=1 \hfil \else \ifodd\pageno
           \hfil {\eightpoint\copy\sectionbox\kern8mm\number\pageno}\else
           {\eightpoint\number\pageno\kern8mm\copy\titlebox}\hfil \fi\fi}
\footline={\hfil}
\headline={\folio}

\def\titlerunning#1{\setbox\titlebox\hbox{\eightpoint #1}}
\def\title#1{\noindent\hfil$\smash{\hbox{\seventeenpointbf #1}}$\hfil
             \titlerunning{#1}\medskip}

\newcount\numbersection \numbersection=-1
\def\sectionrunning#1{\setbox\sectionbox\hbox{\eightpoint #1}}
\def\section#1{%
  \par\vskip0.666cm\penalty -100
  \vbox{\baselineskip=14.4pt\noindent{{\twelvepointbf #1}}}
  \vskip2pt
  \penalty 500
  \advance\numbersection by 1
  \sectionrunning{#1}}

\def\subsection#1|{%
  \par\vskip0.5cm\penalty -100
  \vbox{\noindent{{\tenpointbf #1}}}
  \vskip1pt
  \penalty 500}

\newcount\numberindex \numberindex=0  
\def\index#1#2{%
  \advance\numberindex by 1
  \immediate\write1{\string\def \string\IND #1%
     \romannumeral\numberindex \string{%
     \noexpand#2 \string\dotfill \space \string\S \number\numbersection, 
     p.\string\ \space\number\pageno \string}}}

\newdimen\itemindent \itemindent=\parindent

\def\item#1{\par\noindent\hangindent\itemindent%
            \rlap{#1}\kern\itemindent\ignorespaces}
\def\itemitem#1{\par\noindent\hangindent2\itemindent%
            \kern\itemindent\rlap{#1}\kern\itemindent\ignorespaces}
\def\itemitemitem#1{\par\noindent\hangindent3\itemindent%
            \kern2\itemindent\rlap{#1}\kern\itemindent\ignorespaces}

\long\def\claim#1|#2\endclaim{\par\vskip 5pt\noindent 
{\tenpointbf #1.}\ {\it #2}\par\vskip 5pt}

\def\proof{\noindent{\it Proof}}

\def\today{\ifcase\month\or
January\or February\or March\or April\or May\or June\or July\or August\or
September\or October\or November\or December\fi \space\number\day,
\number\year}

\catcode`\@=11
\newcount\@tempcnta \newcount\@tempcntb 
\def\timeofday{{%
\@tempcnta=\time \divide\@tempcnta by 60 \@tempcntb=\@tempcnta
\multiply\@tempcntb by -60 \advance\@tempcntb by \time
\ifnum\@tempcntb > 9 \number\@tempcnta:\number\@tempcntb
  \else\number\@tempcnta:0\number\@tempcntb\fi}}
\catcode`\@=12

\def\bibitem#1&#2&#3&#4&%
{\hangindent=1.66cm\hangafter=1
\noindent\rlap{\hbox{\eightpointbf #1}}\kern1.66cm{\rm #2}{\it #3}{\rm #4.}} 


\def\bC{{\Bbb C}}

\def\bQ{{\Bbb Q}}
\def\bR{{\Bbb R}}



\def\cO{{\Cal O}}


\def\square{{\hfill \hbox{
\vrule height 1.453ex  width 0.093ex  depth 0ex
\vrule height 1.5ex  width 1.3ex  depth -1.407ex\kern-0.1ex
\vrule height 1.453ex  width 0.093ex  depth 0ex\kern-1.35ex
\vrule height 0.093ex  width 1.3ex  depth 0ex}}}
\def\bigsquare{{\kern-0.3ex\hbox{
\vrule height 1.7ex  width 0.093ex  depth 0ex\kern-0.093ex
\vrule height 1.8ex  width 1.7ex  depth -1.707ex\kern-0.093ex
\vrule height 1.7ex  width 0.093ex  depth 0ex\kern-1.65ex
\vrule height 0.093ex  width 1.6ex  depth 0ex}\kern0.3ex}}
\def\qed{\phantom{$\quad$}\hfill$\square$\medskip}
\def\hexnbr#1{\ifnum#1<10 \number#1\else
 \ifnum#1=10 A\else\ifnum#1=11 B\else\ifnum#1=12 C\else
 \ifnum#1=13 D\else\ifnum#1=14 E\else\ifnum#1=15 F\fi\fi\fi\fi\fi\fi\fi}
\def\msatype{\hexnbr\msafam}
\def\msbtype{\hexnbr\msbfam}
\mathchardef\restriction="3\msatype16   
\mathchardef\compact="3\msatype62
\mathchardef\smallsetminus="2\msbtype72   \let\ssm\smallsetminus
\mathchardef\subsetneq="3\msbtype28
\mathchardef\supsetneq="3\msbtype29
\mathchardef\leqslant="3\msatype36   \let\le\leqslant
\mathchardef\geqslant="3\msatype3E   \let\ge\geqslant
\mathchardef\ltimes="2\msbtype6E
\mathchardef\rtimes="2\msbtype6F


\let\wt=\widetilde
\let\wh=\widehat
\let\text=\hbox
\def\build#1|#2|#3|{\mathrel{\mathop{\null#1}\limits^{#2}_{#3}}}


\def\Vol{\mathop{\rm Vol}\nolimits}

\def\Sing{\mathop{\rm sing}\nolimits}

\def\ddbar{{\partial\overline\partial}}


\def\BC{{\rm BC}}

\def\NS{{\rm NS}}
\def\DNS{{\rm DNS}}


\title{A converse to the}
\smallskip
\title{Andreotti-Grauert theorem}
\titlerunning{A converse to the Andreotti-Grauert theorem}
\medskip
\centerline{\twelvebf Jean-Pierre Demailly}
\medskip
\centerline{Universit\'e de Grenoble I, D\'epartement de Math\'ematiques}
\centerline{Institut Fourier, 38402 Saint-Martin d'H\`eres, France}
\centerline{{\it e-mail\/}: {\tt demailly@fourier.ujf-grenoble.fr}}

\smallskip
\vskip35pt
\line{\hfill \it Dedicated to Professor Nguyen Thanh Van}
\vskip20pt
\noindent
{\bf Abstract.}
The goal of this paper is to show that there are strong relations
between  certain Monge-Amp\`ere integrals appearing in holomorphic Morse
inequa\-lities, and asymptotic cohomology estimates for tensor
powers of holomorphic line bundles. Especially, we prove that these
relations hold without restriction for projective surfaces, and in
the special case of the volume, i.e.\ of asymptotic $0$-cohomology,
for all projective manifolds. These results can be seen as a 
partial converse to the Andreotti-Grauert vanishing theorem.
\medskip

\noindent
{\bf R\'esum\'e.}
Le but de ce  travail est de montrer qu'il y a des relations fortes
entre certaines int\'egrales de Monge-Amp\`ere apparaissant dans les
in\'egalit\'es de Morse holomorphes, et les
estimations asymptotiques de cohomologie pour les fibr\'es holomorphes
en droites. En particulier, nous montrons que ces relations sont
satisfaites sans restriction pour toutes les surfaces projectives,
et dans le cas particulier du volume, c'est-\`a-dire de la
$0$-cohomologie asymptotique, pour toutes les vari\'et\'es projectives.
Ces r\'esultats peuvent \^etre vus comme une r\'eciproque
partielle au th\'eor\`eme d'annulation d'Andreotti-Grauert.

\medskip\noindent
{\bf Mathematics Subject Classification 2010.} 32F07, 14B05, 14C17
\medskip\noindent
{\bf Key words.}
Holomorphic Morse inequalities, Monge-Amp\`ere integrals,
asymptotic cohomology groups, Riemann-Roch formula,
hermitian metrics, Chern curvature tensor, approximate Zariski
decomposition.
\medskip\noindent
{\bf Mots-cl\'es.}
In\'egalit\'es de Morse holomorphes, int\'egrales de Monge-Amp\`ere, 
formule de Riemann-Roch, m\'etriques hermitiennes, courbure de Chern, 
d\'ecompo\-sition de Zariski approch\'ee.
\vfill\eject

\section{1. Main results}

Throughout this paper, $X$ denotes a compact complex manifold, 
\hbox{$n=\dim_\bC X$} its complex dimension and $L\to X$ a holomorphic
line bundle. In order to estimate the growth of cohomology groups, it 
is interesting to consider appropriate ``asymptotic cohomology functions''. 
Following partly notation and concepts introduced by A.~K\"uronya 
[K\"ur06, FKL07], we introduce

\claim (1.1) Definition|
\vskip3pt
{\itemindent=6.5mm
\item{\rm (i)} The $q$-th asymptotic cohomology functional is defined as
$$
\widehat h^q(X,L):=
\limsup_{k\to+\infty}~{n!\over k^n}h^q(X,L^{\otimes k}).
$$
\item{\rm (ii)} The $q$-th asymptotic holomorphic Morse sum of $L$ is
$$
\widehat h^{\,\leq q}(X,L):=
\limsup_{k\to+\infty}~{n!\over k^n}\sum_{0\le j\le q}(-1)^{q-j}h^j(X,L^{\otimes k}).
$$}
\endclaim

When the $\limsup$'s are limits, we have the obvious relation 
$$
\widehat h^{\,\leq q}(X,L)=\sum_{0\le j\le q}(-1)^{q-j}\widehat h^j(X,L).
$$
Clearly, Definition 1.1 can also be given for a $\bQ$-line bundle $L$
or a $\bQ$-divisor $D$, and in the case $q=0$ one gets what is usually
called the volume of $L$, namely
$$
\Vol(X,L)=\widehat h^0(X,L)=
\limsup_{k\to+\infty}~{n!\over k^n}h^0(X,L^{\otimes k}).\leqno(1.2)
$$
(see also [DEL00], [Bou02], [Laz04]). It has been shown in [K\"ur06] for the 
projective case and in [Dem10] in general that the $\smash{\widehat h}^q$
functional induces a continuous map 
$$
\DNS_\bR(X)\ni \alpha\mapsto \widehat h^q(X,\alpha)
$$
defined on the ``divisorial Neron-Severi space'' $\DNS_\bR(X)
\subset H^{1,1}_\BC(X,\bR)$, i.e.\ the vector space spanned by
real linear combinations of classes of divisors in the real Bott-Chern 
cohomology group of bidegree $(1,1)$. Here $H^{p,q}_\BC(X,\bC)$ is defined
as the quotient of $d$-closed $(p,q)$-forms by $\smash{\ddbar}$-exact 
$(p,q)$-forms, and there is a natural conjugation 
$H^{p,q}_\BC(X,\bC)\to H^{q,p}_\BC(X,\bC)$ which allows us to speak
of real classes when $q=p$. The functional
$\smash{\widehat h}^q$ is in fact locally Lipschitz continous on
$\DNS_\bR(X)$, and can be obtained as a limit (not just a limsup) on all
those classes. Notice that $H^{p,q}_\BC(X,\bC)$ coincides with
the usual Dolbeault cohomology group $H^{p,q}(X,\bC)$ when $X$ is K\"ahler,
and that $\DNS_\bR(X)$ coincides with the usual N\'eron-Severi space
$$
\NS_\bR(X)=\bR\otimes_\bQ\big(H^2(X,\bQ)\cap H^{1,1}(X,\bC)\big)
$$
when $X$ is projective.  
It follows from holomorphic Morse inequalities (cf.\ [Dem85], [Dem91])
that asymptotic cohomology can be compared with certain 
Monge-Amp\`ere integrals.

\claim (1.3) Theorem {\rm ([Dem85])}|For every holomorphic line bundle 
$L$ on a compact complex manifold $X$, one has the ``weak Morse inequality''
\vskip3pt
{\itemindent=6.5mm
\item{\rm (i)}\kern50pt
$\displaystyle
\widehat h^q(X,L)\le 
\inf_{u\in c_1(L)}\int_{X(u,q)}(-1)^qu^n$\vskip3pt
\noindent
where $u$ runs over all smooth $d$-closed $(1,1)$-forms which belong to
the cohomology class $c_1(L)\in H^{1,1}_\BC(X,\bR)$, and $X(u,q)$ is the
open set
$$
X(u,q):=\big\{z\in X\,;\;\hbox{$u(z)$ has signature $(n-q,q)$}\big\}.
$$
Moreover, if \hbox{$X(u,\le q):=\bigcup_{0\le j\le q}X(u,j)$}, one
has the ``strong Morse inequality''
\vskip3pt
\item{\rm (ii)}\kern50pt
$\displaystyle
\widehat h^{\,\leq q}(X,L)\le 
\inf_{u\in c_1(L)}\int_{X(u,\le q)}(-1)^qu^n$.\vskip0pt}
\endclaim

It is a natural problem to ask whether the inequalities (1.3)~(i) and
(1.3)~(ii) might not always be equalities. These questions are
strongly related to the Andreotti-Grauert vanishing theorem [AG62]. A
well-known variant of this theorem says that if for some integer $q$
and some $u\in c_1(L)$ the form $u(z)$ has at least $n-q+1$ positive
eigenvalues everywhere (so that \hbox{$X(u,\ge q)= \bigcup_{j\ge
  q}X(u,j)=\emptyset$}), then $H^j(X,L^{\otimes k})=0$ for $j\ge q$
and~$k\gg 1$. We are asking here whether conversely the knowledge that
cohomology groups are asymptotically small in a certain degree $q$ implies the
existence of a hermitian metric on $L$ with suitable curvature, i.e.\
no $q$-index points or only a very small amount of such.

The first goal of this note is to prove that the answer is positive
in the case of the volume functional (i.e.\ in the case of degree $q=0$),
at least when $X$ is projective algebraic.

\claim (1.4) Theorem|Let $L$ be a holomorphic line bundle on a projective
algebraic manifold. then
$$
\Vol(X,L)=\inf_{u\in c_1(L)}\int_{X(u,0)}u^n.
$$
\endclaim

The proof relies mainly on five ingredients: (a) approximate Zariski
decomposition for a K\"ahler current $T\in c_1(L)$ (when $L$ is big), i.e.\
a decomposition $\mu^*T=[E]+\beta$
where $\mu:\smash{\widetilde X}\to X$ is a modification, $E$ an 
exceptional divisor and $\beta$ a K\"ahler metric on $\smash{\widetilde X}\,$; 
(b) the characterization  of the pseudoeffective cone ([BDPP04]), 
and the orthogonality estimate 
$$E\cdot\beta^{n-1}\le C\big({\rm Vol}(X,L)-\beta^n\big)^{1/2}$$
proved as an intermediate step of that characterization; 
(c) properties of solutions of Laplace equations to get smooth
approximations of $[E]\,$; (d) log concavity of the Monge-Amp\`ere 
operator$\,$; and finally (e) birational invariance of the Morse infimums. 
In the case of higher cohomology groups, we have been able to treat 
only the case of projective surfaces~:

\claim (1.5) Theorem|Let $L\to X$ be a holomorphic line bundle on
a complex projective surface. Then both weak and strong 
inequalities {\rm(1.3)~(i)} and {\rm(1.3)~(ii)} 
are equalities for $q=0,\,1,\,2$, and the $\limsup$'s involved in
$\smash{\widehat h}^q(X,L)$ and $\smash{\widehat h}^{\,\leq q}(X,L)$
are limits.
\endclaim

Thanks to the Serre duality and the Riemann-Roch formula, the
(in)equality for a given $q$ is equivalent to the (in)equality for
$n-q$. Therefore, on surfaces, the only substantial case which still
has to be checked in addition to Theorem 1.4 is the case $q=1\,$: this
is done by using Grauert' criterion that the intersection matrix
$(E_i\cdot E_j)$ is negative definite for every exceptional divisor
$E=\sum c_jE_j$. Our statements are of course trivial on curves since
the curvature of any holomorphic line bundle can be taken to be
constant with respect to any given hermitian metric.

\claim (1.6) Remark|{\rm It is interesting put these results in perspective
with the algebraic version of holomorphic Morse inequalities proved in 
[Dem94] (see  also [Siu93] and [Tra95] for related ideas, and [Ang94] for an 
algebraic proof).  When $X$ is projective, the algebraic Morse inequalities
used in combination with the birational invariance of the Morse integrals
(cf.\ section~2) imply the inequalities
{\itemindent=6.5mm
\vskip5pt
\item{\rm (i)}\kern20pt
$\displaystyle\inf_{u\in c_1(L)}\int_{X(u,q)}(-1)^qu^n
\leq \inf_{\mu^*(L)\simeq\cO(A-B)}{n\choose q}A^{n-q}B^q\,,$
\vskip5pt
\item{\rm (ii)}\kern20pt
$\displaystyle\inf_{u\in c_1(L)}\int_{X(u,\le q)}(-1)^qu^n
\leq \inf_{\mu^*(L)\simeq\cO(A-B)}\sum_{0\le j\le q}(-1)^{q-j}
{n\choose j}A^{n-j}B^j\,,$\vskip5pt}
\noindent
where the infimums on the right hand side are taken over all
modifications $\mu:\smash{\widetilde X}\to X$ and all decompositions
$\mu^*L=\cO(A-B)$ of $\mu^*L$ as a difference of two nef $\bQ$-divisors 
$A,\,B$ on $\smash{\widetilde X}$. In case $A$ and $B$ are ample, the
proof simply consists of taking positive curvature forms
$\Theta_{\cO(A),h_A}$, $\Theta_{\cO(B),h_B}$ on $\cO(A)$ and $\cO(B)$,
and evaluating the Morse integrals with $u=\Theta_{\cO(A),h_A}-
\Theta_{\cO(B),h_B}\,$; the 
general case follows by approximating the nef divisors $A$ and $B$ by 
ample divisors $A+\varepsilon H$ and $B+\varepsilon H$ with $H$ ample and
$\varepsilon>0$, see [Dem94]. Again, a natural question is to know 
whether these infimums derived from algebraic intersection numbers
are equal to the asymptotic cohomology functionals 
$\smash{\widehat h}^q(X,L)$ and $\smash{\widehat h}^{\,\leq q}(X,L)$.
A positive answer would of course automatically yield a positive answer 
to the equality cases in 1.3~(i) and 1.3~(ii). However, the Zariski 
decompositions involved in our proofs of the ``analytic equality case'' 
produces certain effective exceptional divisors which are not nef. It
is unclear how to write those effective divisors as a difference of nef 
divisors. This fact raises some doubts upon the sufficiency of taking
merely differences of nef divisors in the infimums 1.6~(i) and 1.6~(ii).\qed}
\endclaim

I warmly thank Burt Totaro for stimulating discussions in connection with
his recent work [Tot10]. 

\section{2. Invariance by modification}

It is easy to check that the asymptotic cohomology function is invariant
by modification, namely that for every modification $\mu:\wt X\to X$ and
every line bundle $L$ we have
$$
\wh h^q(X,L)=\wh h^q(\wt X,\mu^*L).\leqno(2.1)
$$
In fact the Leray's spectral sequence provides an $E_2$ term
$$
E_2^{p,q}=H^p(X,R^q\mu_*\cO_{\wt X}(\mu^*L^{\otimes k}))=
H^p(X,\cO_X(L^{\otimes k})\otimes R^q\mu_*\cO_{\wt X}).
$$
Since $R^q\mu_*\cO_{\wt X}$ is equal to $\cO_X$ for $q=0$ and is supported
on a proper analytic subset of $X$ for $q\ge 1$, one infers that
$h^p(X,\cO_X(L^{\otimes k}\otimes R^q\mu_*\cO_{\wt X})=O(k^{n-1})$ for all
$q\ge 1$. The spectral sequence implies that
$$
h^q(X,L^{\otimes k})-\wh h^q(\wt X,\mu^*L^{\otimes k})=O(k^{n-1}).
$$
We claim that the Morse integral infimums are also invariant by modification.

\claim (2.2) Proposition|Let $(X,\omega)$ be a compact K\"ahler manifold,
$\alpha\in H^{1,1}(X,\bR)$ a real cohomology class and 
$\mu:\wt X\to X$ a modification. Then
$$
\leqalignno{
&\strut\kern20pt\inf_{u\in \alpha}\int_{X(u,q)}(-1)^qu^n
=\inf_{v\in \mu^*\alpha}\int_{X(v,q)}(-1)^qv^n,
&\hbox{\rm(i)}\cr
&\strut\kern20pt\inf_{u\in \alpha}\int_{X(u,\le q)}(-1)^qu^n
=\inf_{v\in \mu^*\alpha}\int_{X(v,\le q)}(-1)^qv^n.
&\hbox{\rm(ii)}\cr}
$$
\endclaim

\proof. Given $u\in\alpha$ on $X$, we obtain Morse integrals with the same 
values by taking $v=\mu^*u$ on~$\wt X$, hence the infimum on $\wt X$ 
is smaller or equal to what is on~$X$. Conversely, we have to show 
that given a smooth representative 
$v\in\mu^*\alpha$ on $\smash{\wt X}$, one can find a smooth representative 
$u\in X$ such that the Morse integrals do not differ much. We can
always assume that $\wt X$ itself is K\"ahler, since by Hironaka [Hir64]
any modification $\smash{\wt X}$ is dominated by a composition of blow-ups
of $X$. Let us fix some $u_0\in\alpha$ and write
$$
v=\mu^*u_0+dd^c\varphi
$$
where $\varphi$ is a smooth function on $\wt X$. We adjust $\varphi$ by a
constant in such a way that $\varphi\ge 1$ on $\wt X$. There exists an analytic
set $S\subset X$ such that $\mu:\wt X\ssm\mu^{-1}(S)\to X\ssm S$ is a
biholomorphism, and a quasi-psh function $\psi_S$ which is smooth on
$X\ssm S$ and has $-\infty$ logarithmic poles on $S$ (see e.g.\
[Dem82]). We define
$$
\wt u=
\mu^*u_0+dd^c\max\nolimits_{\varepsilon_0}(\varphi+\delta\,\psi_S\circ\mu~,~0)
=v+dd^c\max\nolimits_{\varepsilon_0}(\delta\,\psi_S\circ\mu~,~-\varphi)
\leqno(2.3)
$$
where $\max\nolimits_{\varepsilon_0}$, $0<\varepsilon_0<1$, is a 
regularized max function and 
$\delta>0$ is very small. By construction $\wt u$ coincides with $\mu^*u_0$ in
a neighborhood of $\mu^{-1}(S)$ and therefore $\wt u$ descends to a smooth 
closed $(1,1)$-form $u$ on $X$ which coincides with $u_0$ near~$S$, so that 
$\wt u=\mu^* u$. Clearly $\wt u$ converges uniformly to $v$ on every compact 
subset of $\wt X\ssm \mu^{-1}(S)$ as $\delta\to 0$, so we only
have to show that the Morse integrals are small (uniformly in $\delta$) 
when restricted to a
suitable small neighborhood of the exceptional set $E=\mu^{-1}(S)$.
Take a sufficiently large K\"ahler metric $\wt\omega$ on $\wt X$ such that
$$
-{1\over 2}\wt\omega\le v\le{1\over 2}\wt\omega,\quad
-{1\over 2}\wt\omega\le dd^c\varphi\le{1\over 2}\wt\omega,\quad
-\wt\omega\le dd^c\psi_S\circ\mu.
$$
Then $\wt u\ge -\wt\omega$ and 
$\wt u\le \wt\omega+\delta\,dd^c\psi_S\circ\mu$ everywhere on $\wt X$.
As a consequence
$$
|\wt u^n|\le \big(\wt\omega+\delta(\wt\omega+dd^c\psi_S\circ\mu)\big)^n
\le\wt\omega^n+n\delta(\wt\omega+dd^c\psi_S\circ\mu)\wedge
\big(\wt\omega+\delta(\wt\omega+dd^c\psi_S\circ\mu)\big)^{n-1}
$$
thanks to the inequality $(a+b)^n\le a^n+nb(a+b)^{n-1}$. For any neighborhood
$V$ of $\mu^{-1}(S)$ this implies
$$
\int_V|\wt u^n|\le \int_V\wt\omega^n+n\delta(1+\delta)^{n-1}\int_{\wt X}
\wt\omega^n
$$
by Stokes formula. We thus see that the integrals are small if $V$ and
$\delta$ are small. The reader may be concerned that Monge-Amp\`ere
integrals were used with an unbounded potential $\psi_S$, but in fact,
for any given $\delta$, all the above formulas and estimates are still
valid when we replace $\psi_S$ by 
$\max_{\varepsilon_0}(\psi_S, -(M+2)/\delta)$ with
$M=\max_{\wt X}\varphi$, especially formula (2.3) shows that
the form $\wt u$ is unchanged.
Therefore our calculations can be handled by using
merely smooth potentials.\qed

\section{3. Proof of the infimum formula for the volume}

We have to show that
$$
\inf_{u\in c_1(L)}\int_{X(u,0)}u^n\le \Vol(X,L)\leqno(3.1)
$$
Let us first assume that $L$ is a big line bundle, i.e.\ 
that $\Vol(X,L)>0$. Then it is known by [Bou02] that
$\Vol(X,L)$ is obtained as the supremum
of $\smash{\int_{X\ssm \Sing(T)}T^n}$ for K\"ahler currents 
$T=-{i\over 2\pi}\ddbar h$ with analytic
singularities in $c_1(L)$; this means that locally $h=e^{-\varphi}$
where $\varphi$ is a strictly plurisubharmonic function which has the
same singularities as $c\log\sum|g_j|^2$ where $c>0$ and the
$g_j$ are holomorphic functions. By [Dem92], there exists a
blow-up $\mu:\widetilde X\to X$ such that $\mu^*T=[E]+\beta$ where $E$ is a
normal crossing divisor on $\widetilde X$ and $\beta\ge 0$
smooth. Moreover, by [BDPP04] we have the orthogonality estimate
$$
[E]\cdot \beta^{n-1}=\int_{E}\beta^{n-1}\le 
C\big(\Vol(X,L)-\beta^n\big)^{1/2},\leqno(3.2)
$$
while
$$
\beta^n=\int_{\wt X}\beta^n=\int_{X\ssm \Sing(T)}T^n\quad
\hbox{approaches $\Vol(X,L)$}.\leqno(3.3)
$$
In other words, $E$ and $\beta$ become ``more and more orthogonal'' as
$\beta^n$ approaches the volume (these properties are summarized
by saying that $\mu^*T=[E]+\beta$ defines an approximate Zariski 
decomposition of $c_1(L)$, cf.\ also [Fuj94]). By subtracting to $\beta$
a small linear combination of the exceptional divisors and increasing
accordingly the coefficients of $E$, we can achieve that the
cohomology class $\{\beta\}$ contains a positive definite form
$\beta'$ on $\wt X$ (i.e.\ the fundamental 
form of a K\"ahler metric); we refer e.g.\ to ([DP04], proof of
Lemma~3.5) for details. This means that we can replace $T$ by
a cohomologous current such that the corresponding form $\beta$ is 
actually a K\"ahler metric, and we will assume for simplicity of notation
that this situation occurs right away for~$T$. Under this assumption,
there exists a smooth closed $(1,1)$-form $v$ belonging to the 
Bott-Chern cohomology class of $[E]$, such that we have identically
$(v-\delta\beta)\wedge\beta^{n-1}=0$ where
$$
\delta={[E]\cdot\beta^{n-1}\over \beta^n}\le C'(\Vol(X,L)-\beta^n\big)^{1/2}
\leqno(3.4)
$$
for some constant $C'>0$. In fact, given an arbitrary smooth
representative $v_0\in\{[E]\}$, the existence of $v=v_0+i\ddbar\psi$ amounts 
to solving a Laplace equation $\Delta\psi=f$ with respect to the K\"ahler 
metric $\beta$, and the choice of $\delta$ ensures that we have 
$\int_X f\,\beta^n=0$ and hence that the equation is solvable.
Then $\wt u:=v+\beta$ is a smooth closed $(1,1)$-form in the cohomology
class $\mu^* c_1(L)$, and its eigenvalues with respect to $\beta$ are
of the form $1+\lambda_j$ where $\lambda_j$ are the eigenvalues of $v$. The
Laplace equation is equivalent to the identity $\sum\lambda_j=n\delta$. 
Therefore
$$
\sum_{1\le j\le n}\lambda_j\le C''(\Vol(X,L)-\beta^n\big)^{1/2}.
\leqno(3.5)
$$
The inequality between arithmetic means and geometry means implies
$$
\prod_{1\le j\le n}(1+\lambda_j)\le 
\Big(1+{1\over n}\sum_{1\le j\le n}\lambda_j\Big)^n
\le 1+C_3(\Vol(X,L)-\beta^n\big)^{1/2}
$$
whenever all factors $(1+\lambda_j)$ are nonnegative. By 2.2~(i) we get
$$
\eqalign{
\inf_{u\in c_1(L)}\int_{X(u,0)}u^n&\le\int_{\wt X(\wt u,0)}\wt u^n\cr
&\le 
\int_{\wt X}\beta^n\big(1+C_3(\Vol(X,L)-\beta^n\big)^{1/2}\big)\cr
&\le \Vol(X,L)+C_4(\Vol(X,L)-\beta^n\big)^{1/2}.\cr}
$$
As $\beta^n$ approches $\Vol(X,L)$, this implies inequality (3.1).
\medskip

We still have to treat the case when $L$ is not big, i.e.\ $\Vol(X,L)=0$. 
Let $A$ be an ample
line bundle and let $t_0\ge 0$ be the infimum of real numbers such that
$L+tA$ is a big $\bQ$-line bundle for $t$ rational, $t>t_0$. The
continuity of the volume function implies that $0<\Vol(X,L+tA)\le\varepsilon$ 
for $t>t_0$ sufficiently close to~$t_0$. By what we have just proved, there
exists a smooth form $u_t\in c_1(L+tA)$ such that $\int_{X(u_t,0)}u_t^n\le
2\varepsilon$. Take a K\"ahler metric $\omega\in c_1(A)$ and define
$u=u_t-t\omega$. Then clearly
$$
\int_{X(u,0)}u^n\le \int_{X(u_t,0)}u_t^n\le 2\varepsilon,
$$
hence
$$
\inf_{u\in c_1(L)}\int_{X(u,0)}u^n=0.
$$
Inequality (3.1) is now proved in all cases.\qed

\section{4. Estimate of the first cohomology group on a projective surface}

Assume first that $L$ is a big line bundle on a projective non singular
variety~$X$. As in section 3, we can find an approximate Zariski decomposition,
i.e.\ a blow-up $\mu:\wt X\to X$ and a current $T\in c_1(L)$ such
$\mu^*T=[E]+\beta$, where $E$ an effective divisor and $\beta$ a K\"ahler 
metric on $\wt X$ such that
$$
\Vol(X,L)-\eta<\beta^n<\Vol(X,L),\qquad\eta\ll 1.\leqno(4.1)
$$
(On a projective surface, one can even get exact Zariski decomposition, 
but we want to remain general at this point, so that the arguments might
possibly be applied later for arbitrary dimension).
By blowing-up further, we may even assume that $E$ is a normal crossing 
divisor. We select a hermitian metric $h$ on $\cO(E)$ and take
$$
u_\varepsilon={i\over 2\pi}\ddbar\log(|\sigma_E|_h^2+\varepsilon^2)+
\Theta_{\cO(E),h}+\beta~~\in~~\mu^*c_1(L)
\leqno(4.2)
$$
where $\sigma_E\in H^0(\widetilde X,\cO(E))$ is the canonical section
and $\Theta_{\cO(E),h}$ the Chern curvature form. Clearly, by the
Lelong-Poincar\'e equation, $u_\varepsilon$ converges to $[E]+\beta$ in 
the weak topology as $\varepsilon\to 0$. Straightforward calculations yield
$$
u_\varepsilon={i\over 2\pi}{\varepsilon^2D^{1,0}_h\sigma_E\wedge 
\overline{D^{1,0}_h\sigma_E}\over(\varepsilon^2+|\sigma_E|^2)^2}
+{\varepsilon^2\over\varepsilon^2+|\sigma_E|^2}\Theta_{E,h}+\beta.
$$
The first term converges to $[E]$ in the weak topology, while the second,
which is close to $\Theta_{E,h}$ near~$E$, converges pointwise everywhere to $0$
on $\smash{\widetilde X}\ssm E$. A simple asymptotic analysis shows that 
$$
\Big({i\over 2\pi}{\varepsilon^2D^{1,0}_h\sigma_E\wedge 
\overline{D^{1,0}_h\sigma_E}\over(\varepsilon^2+|\sigma_E|^2)^2}
+{\varepsilon^2\over\varepsilon^2+|\sigma_E|^2}\Theta_{E,h}\Big)^p\to
[E]\wedge\Theta_{E,h}^{p-1}
$$
in the weak topology for $p\ge 1$, hence
$$
\lim_{\varepsilon\to 0}u_\varepsilon^n=
\beta^n+\sum_{p=1}^n{n\choose p}[E]\wedge \Theta_{E,h}^{p-1}
\wedge\beta^{n-p}.\leqno(4.3)
$$
In arbitrary dimension, the signature of $u_\varepsilon$ is hard to evaluate,
and it is also non trivial to decide the sign of the limiting measure
$\lim u_\varepsilon^n$. However, when $n=2$, we get the simpler formula
$$
\lim_{\varepsilon\to 0}u_\varepsilon^2=\beta^2+2[E]\wedge\beta
+[E]\wedge\Theta_{E,h}.
$$
In this case, $E$ can be assumed to be an exceptional divisor (otherwise
some part of it would be nef and could be removed from the poles of $T$).
Hence the matrix $(E_j\cdot E_k)$ is negative definite and we can find 
a smooth hermitian
metric $h$ on $\cO(E)$ such that $(\Theta_{E,h})_{|E}<0$, i.e.\
$\Theta_{E,h}$ has one negative eigenvalue everywhere along~$E$. 

\claim (4.4) Lemma|One can adjust the metric $h$ of $\cO(E)$ in such a
way that $\Theta_{E,h}$ is negative definite 
on a neighborhood of the support $|E|$ of the exceptional divisor,
and $\Theta_{E,h}+\beta$ has signature $(1,1)$ there.
$($We do not care about the signature far away from $|E|)$.
\endclaim

\proof. At a given point $x_0\in X$, let us fix coordinates and a positive
quadratic form $q$ on $\bC^2$. If we put
$\psi_\varepsilon(z)=\varepsilon\chi(z)\log(1+\varepsilon^{-1}q(z))$
with a suitable cut-off function~$\chi$,
then the Hessian form of $\psi_\varepsilon$ is equal to $q$ at $x_0$ and
decays rapidly to $O(\varepsilon\log\varepsilon)|dz|^2$ away from $x_0$. In 
this way, after multiplying $h$ with $e^{\pm \psi_\varepsilon(z)}$, we can replace
the curvature $\Theta_{E,h}(x_0)$ with $\Theta_{E,h}(x_0)\pm q$ without
substantially modifying the form away from~$x_0$. This allows to adjust
$\Theta_{E,h}$ to be equal to (say) $-{1\over 4}\beta(x_0)$ at any singular
point $x_0\in E_j\cap E_k$ in the support of $|E|$, while keeping
$\Theta_{E,h}$ negative definite along~$E$. In order to adjust the
curvature at smooth points $x\in|E|$, we replace the metric $h$ with
$h'(z)= h(z)\exp(-c(z)|\sigma_E(z)|^2)$. Then the curvature form 
$\Theta_{E,h}$ is
replaced by $\Theta_{E,h'}(x)=\Theta_{E_h}(x)+c(x)|d\sigma_E|^2$ at 
$x\in|E|$ (notice that $d\sigma_E(x)=0$ if $x\in\Sing|E|$), and we can
always select a real function $c$ so that 
$\Theta_{E,h'}$ is negative definite with one negative
eigenvalue between $-1/2$ and~$0$ at any point of $|E|$. Then 
$\Theta_{E,h'}+\beta$ has signature $(1,1)$ near~$|E|$.\qed

With this choice of the metric, we see that for $\varepsilon>0$ small,
the sum
$$
{\varepsilon^2\over\varepsilon^2+|\sigma_E|^2}\Theta_{E,h}+\beta
$$
is of signature $(2,0)$ or $(1,1)$ (or degenerate of signature $(1,0)$),
the non positive definite points being concentrated in a neighborhood
of~$E$. In particular the index set $X(u_\varepsilon,2)$ is empty, and 
also
$$
u_\varepsilon\le
{i\over 2\pi}{\varepsilon^2D^{1,0}_h\sigma_E\wedge 
\overline{D^{1,0}_h\sigma_E}\over(\varepsilon^2+|\sigma_E|^2)^2}+\beta
$$
on a neighborhood $V$ of $|E|$, while $u_\varepsilon$ converges uniformly
to $\beta$ on $\wt X\ssm V$. This implies that
$$
\beta^2\le\liminf_{\varepsilon\to 0}\int_{X(u_\varepsilon,0)}u_\varepsilon^2
\le\limsup_{\varepsilon\to 0}\int_{X(u_\varepsilon,0)}u_\varepsilon^2
\le \beta^2+2\beta\cdot E.
$$
Since $\int_{\wt X}u_\varepsilon^2=L^2=\beta^2+2\beta\cdot E+E^2$ we conclude
by taking the difference that
$$
-E^2-2\beta\cdot E
\le\liminf_{\varepsilon\to 0}\int_{X(u_\varepsilon,1)}-u_\varepsilon^2
\le\limsup_{\varepsilon\to 0}\int_{X(u_\varepsilon,1)}-u_\varepsilon^2
\le -E^2.
$$
Let us recall that $\beta\cdot E\le C(\Vol(X,L)-\beta^2)^{1/2}=0(\eta^{1/2})$ 
is small by (4.1) and the orthogonality estimate. The asymptotic 
cohomology is given here
by $\wh h^2(X,L)=0$ since $h^2(X,L^{\otimes k})=H^0(X,K_X\otimes L^{\otimes -k})=
0$ for $k\ge k_0$, and we have by Riemann-Roch
$$
\wh h^1(X,L)=\wh h^0(X,L)-L^2=\Vol(X,L)-L^2=-E^2-\beta\cdot E+O(\eta).
$$
Here we use the fact that ${n!\over k^n}h^0(X,L^{\otimes k})$ converges to 
the volume when $L$ is big. All this shows that equality occurs in the Morse 
inequalities~(1.3) when we pass to the infimum. By taking limits in
the Neron-Severi space $\NS_\bR(X)\subset H^{1,1}(X,\bR)$, we
further see that equality occurs as soon as $L$ is pseudo-effective, and the 
same is true if $-L$ is pseudo-effective by Serre duality.
\medskip
It remains to treat the case when neither $L$ nor $-L$ are pseudo-effective.
Then $\wh h^0(X,L)=\wh h^2(X,L)=0$, and asymptotic cohomology appears
only in degree~$1$, with $\wh h^1(X,L)=-L^2$ by Riemann-Roch.
Fix an ample line bundle $A$ and let $t_0>0$ be the infimum of real
numbers such that $L+tA$ is big  for $t$~rational, $t>t_0$, 
resp.\ let $t'_0>0$ be the infimum of real
numbers $t'$ such that $-L+t'A$ is big  for $t'>t'_0$. Then for
$t>t_0$ and $t'>t'_0$, we can find a modification 
$\mu:\wt X\to X$ and currents $T\in c_1(L+tA)$,
$T'\in c_1(-L+t'A)$ such that
$$
\mu^*T=[E]+\beta,\qquad \mu^*T'=[F]+\gamma
$$
where $\beta$, $\gamma$ are K\"ahler forms and $E$, $F$ normal crossing
divisors. By taking a suitable linear combination $t'(L+tA)-t(-L+t'A)$
the ample divisor $A$ disappears, and we get
$$
{1\over t+t'}\Big(t'[E]+t'\beta-t[F]-t\gamma\Big)\in \mu^*c_1(L).
$$
After replacing $E$, $F$, $\beta$, $\gamma$ by suitable multiples, we
obtain an equality
$$
[E]-[F]+\beta-\gamma\in \mu^*c_1(L).
$$
We may further assume by subtracting that the divisors $E$, $F$ have no 
common components.
The  construction shows that $\beta^2\le\Vol(X,L+tA)$ can be taken 
arbitrarily small (as well of course as $\gamma^2$), and the orthogonality
estimate implies that we can assume $\beta\cdot E$ and $\gamma\cdot F$
to be arbitrarily small. Let us introduce metrics $h_E$ on 
$\cO(E)$ and $h_F$ on $\cO(F)$ as in Lemma~4.4, and consider the forms
$$
\eqalign{
u_\varepsilon=
&+{i\over 2\pi}{\varepsilon^2D^{1,0}_{h_E}\sigma_E\wedge 
\overline{D^{1,0}_{h_E}\sigma_E}\over(\varepsilon^2+|\sigma_E|^2)^2}
+{\varepsilon^2\over\varepsilon^2+|\sigma_E|^2}\Theta_{E,h_E}+\beta\cr
&-{i\over 2\pi}{\varepsilon^2D^{1,0}_{h_F}\sigma_F\wedge 
\overline{D^{1,0}_{h_F}\sigma_F}\over(\varepsilon^2+|\sigma_F|^2)^2}
-{\varepsilon^2\over\varepsilon^2+|\sigma_F|^2}\Theta_{F,h_F}-\gamma
~~\in~~\mu^*c_1(L).\cr}
$$
Observe that $u_\varepsilon$ converges uniformly to $\beta-\gamma$ outside
of every neighborhood of $|E|\cup|F|$. Assume that
$\Theta_{E,h_E}<0$ on $V_E=\{|\sigma_E|<\varepsilon_0\}$ and
$\Theta_{F,h_F}<0$ on $V_F=\{|\sigma_F|<\varepsilon_0\}$. On 
$V_E\cup V_F$ we have
$$
u_\varepsilon\le 
{i\over 2\pi}{\varepsilon^2D^{1,0}_{h_E}\sigma_E\wedge 
\overline{D^{1,0}_{h_E}\sigma_E}\over(\varepsilon^2+|\sigma_E|^2)^2}
-{\varepsilon^2\over\varepsilon^2+|\sigma_F|^2}\Theta_{F,h_F}+\beta
+{\varepsilon^2\over\varepsilon_0^2}\Theta_{E,h_E}^+
$$
where $\Theta_{E,h_E}^+$is the positive part of $\Theta_{E,h_E}$ with respect
to $\beta$. One sees immediately that this term is negligible. The first
term is the only one which is not uniformly bounded, and actually it
converges weakly to the current $[E]$. By squaring, we find
$$
\limsup_{\varepsilon\to 0}\int_{X(u_\varepsilon,0)}u_\varepsilon^2\le
\int_{X(\beta-\gamma,0)}(\beta-\gamma)^2+2\beta\cdot E.
$$
Notice that the term $-\smash{\varepsilon^2\over\varepsilon^2+|\sigma_F|^2}\,
\Theta_{F,h_F}$ does not contribute to the limit as
it converges boun\-dedly almost everywhere to $0$, the exceptions
being points of~$|F|$, but this set is
of measure zero with respect to the current $[E]$. Clearly we have
$\int_{X(\beta-\gamma,0)}(\beta-\gamma)^2\le \beta^2$ and therefore
$$
\limsup_{\varepsilon\to 0}\int_{X(u_\varepsilon,0)}u_\varepsilon^2\le
\beta^2+2\beta\cdot E.
$$
Similarly, by looking at $-u_\varepsilon$, we find
$$
\limsup_{\varepsilon\to 0}\int_{X(u_\varepsilon,2)}u_\varepsilon^2\le
\gamma^2+2\gamma\cdot F.
$$
These $\limsup$'s are small and we conclude that the essential part of
the mass is concentrated on the $1$-index set, as desired.\qed

\section{References}
\medskip

{\eightpoint

\bibitem[AG62]&Andreotti, A., Grauert, H.:& Th\'eor\`emes de finitude
  pour la cohomologie des espaces complexes;& Bull.\ Soc.\ Math.\
  France {\bf 90} (1962) 193--259&

\bibitem[Ang95]&Angelini, F.:& An algebraic version of Demailly's
  asymptotic Morse inequalities;& arXiv: alg-geom/9503005,
  Proc.\ Amer.\ Math.\ Soc.\ {\bf 124} (1996) 3265--3269&

\bibitem[Bou02]&Boucksom, S.:& On the volume of a line bundle;&
  Internat.\ J.\ Math.\ {\bf 13} (2002), 1043--1063&

\bibitem[BDPP04]&Boucksom, S., Demailly, J.-P., P\u{a}un, M.,
  Peternell, Th.:& The pseudo-effective cone of a compact K\"ahler
  manifold and varieties of negative Kodaira dimension;& arXiv: 
  math.AG/0405285, see also Proceedings of the ICM 2006 in Madrid&

\bibitem[Dem82]&Demailly, J.-P.:& Estimations $L^2$ pour l'op\'erateur
  $\overline \partial $ d'un fibr\'e vectoriel holomorphe semi-positif
  au dessus d'une vari\'et\'e k\"ahl\'erienne compl\`ete;& Ann.\ Sci.\
  \'Ecole Norm.\ Sup.\ {\bf 15} (1982) 457--511&

\bibitem[Dem85]&Demailly, J.-P.:& Champs magn\'etiques et
  in\'egalit\'es de Morse pour la $d''$-cohomo\-logie;& Ann.\ Inst.\
  Fourier (Grenoble) {\bf 35} (1985), 189--229&

\bibitem[Dem91]&Demailly, J.-P.:& Holomorphic Morse inequalities;&
  Lectures given at the AMS Summer Institute on Complex Analysis held
  in Santa Cruz, July 1989, Proceedings of Symposia in Pure
  Mathematics, Vol.~{\bf 52}, Part~2 (1991), 93--114&

\bibitem[Dem94]&Demailly, J.-P.:& $L^2$ vanishing theorems for positive 
  line bundles and adjunction theory;& arXiv: alg-geom/9410022$\,$;
  Lecture Notes of the CIME Session ``Transcendental methods in
  Algebraic Geometry'', Cetraro, Italy, July 1994, Ed.\
  F.~Catanese, C.~Cili\-berto, Lecture Notes in Math., Vol.~1646, 1--9&

\bibitem[Dem10]&Demailly, J.-P.:& Holomorphic Morse inequalities and
  asymptotic cohomology groups: a tribute to Bernhard Riemann;&
  Milan Journal of Mathematics {\bf 78} (2010) 265--277&

\bibitem[Dem92]&Demailly, J.-P.:& Regularization of closed positive currents 
  and Intersection Theory;& J.\ Alg.\ Geom.\ {\bf 1} (1992), 361--409&

\bibitem[DEL00]&Demailly, J.-P., Ein, L., {\rm and} Lazarsfeld, R.:& A
  subadditivity property of multiplier ideals;& Michigan Math. J. 48
  (2000), 137\--156&

\bibitem[DP04]&Demailly, J.-P., P\u{a}un, M:& Numerical
  characterization of the K\"ahler cone of a compact K\"ahler
  manifold;& arXiv: math.AG/0105176$\,$; Annals of Math.\ {\bf 159} (2004)
  1247--1274&

\bibitem[FKL07]&de Fernex, T., K\"uronya, A., Lazarsfeld, R.:& Higher
  cohomology of divisors on a projective variety;& Math.\ Ann.\ {\bf
  337} (2007) 443--455&

\bibitem[Fuj94]&Fujita, T.: & Approximating Zariski decomposition of
  big line bundles;& Kodai Math.\ J.\ {\bf 17} (1994) 1-\-3&

\bibitem[Hir64]&Hironaka, H.:& Resolution of singularities of an
  algebraic variety over a field of characteristic zero;& Ann.\ of
  Math.\ {\bf 79} (1964) 109--326&

\bibitem[K\"ur06]&K\"uronya, A.:& Asymptotic cohomological functions on
  projective varieties;& Amer.\ J.\ Math.\ {\bf 128} (2006) 1475--1519&

\bibitem[Laz04]&Lazarsfeld, R.:& Positivity in Algebraic Geometry
  I.-II;& Ergebnisse der Mathematik und ihrer Grenzgebiete,
  Vols.\ 48-49, Springer Verlag, Berlin, 2004&

\bibitem[Siu93]&Siu, Y.T.:& An effective Matsusaka big theorem;& Ann.\ Inst.\
  Fourier {\bf 43} (1993) 1387--1405&

\bibitem[Tra95]&Trapani, S.:& Numerical criteria for the positivity of the
  difference of ample divisors;& Math.\ Zeitschrift {\bf 219} (1995) 387--401&

\bibitem[Tot10]&Totaro, B.:& Line bundles with partially vanishing cohomology;&
  July 2010, arXiv: math.AG/1007.3955&

}
\vskip20pt
\noindent
(version of December 17, 2010, printed on \today)
\end